%% file: Zeros_linear_twists_L_functions.tex
\documentclass[10pt]{amsart}

\usepackage{graphicx}
\usepackage[plainpages = false, pdfpagelabels, 
						 bookmarks,
						 bookmarksopen = true,
						 bookmarksnumbered = true,
						 breaklinks = true,
						 linktocpage,
						 colorlinks = true,
						 linkcolor = blue,
						 urlcolor  = blue,
						 citecolor = red,
						 anchorcolor = green,
						 hyperindex = true,
						 hyperfigures
						 ]{hyperref}
\usepackage[utf8]{inputenc}
\usepackage[english]{babel}
\usepackage[T1]{fontenc}

\usepackage{type1cm}

\usepackage{amsmath, amsfonts, amsthm, amscd, amssymb}
\usepackage[shortlabels]{enumitem}
\usepackage[normalem]{ulem}
\usepackage{extarrows}
\usepackage[all]{xy}

\usepackage[a4paper]{geometry}
\usepackage{vmargin}
    \setpapersize{A4}
    \setmarginsrb{26mm}{28mm}{26mm}{25mm}	
    				{8mm}{6mm}{10mm}{12mm}				

\usepackage{abbreviazioni}
\usepackage{stile_en}

\newcommand{\prima}[1]{}

\title{Zeros of linear twists of $L$-functions outside the critical strip}
\author{Mattia Righetti}
\date{}
\address{Dipartimento di Matematica, Universit\`{a} di Genova}
\email{righetti@dima.unige.it}           

\begin{document}

\begin{abstract}
In this note we investigate the existence of zeros of linear twists of $L$-functions outside of the critical strip. In particular, we show that the Lerch zeta function $L(\lambda,\alpha,s)$ has infinitely many zeros for $1<\sigma<1+\eta$, for any $\eta>0$, when $\lambda$ is irrational and $\alpha$ is rational. This settles the question on the existence of zeros of the Lerch zeta functions for $\sigma>1$.\\

\end{abstract}

\maketitle

\input{intro.tex}

\input{exponential_sums.tex}

\input{existence_zeros.tex}
\input{continuity.tex}

\medskip
\begin{acknowledgements}
I would like to thank Professor Alberto Perelli for the many valuable discussions and suggestions, and Professor Giuseppe Molteni for interesting discussions. Furthermore, this research has been conducted with the financial support of a research scholarship from the Department of Mathematics of the University of Genova, which I would like to thank.
\end{acknowledgements} 

\bibliographystyle{amsplain}
\bibliography{biblio}

\end{document}

%% file: intro.tex
\section{Introduction}

In \cite{lerch}, Lerch introduced and studied the function
$$L(\lambda ,\alpha,s)=\sum_{n=0}^\infty \frac{\e{\lambda n}}{(n+\alpha)^s}, \qquad \sigma>1,\quad 0<\alpha,\lambda \leq 1, \quad \e{z}=\exp{2\pi i z}.$$
When $\lambda =1$, the Lerch zeta function reduces to the Hurwitz zeta function $\zeta(s,a)$, so it is well known that it has infinitely many zeros for  for $1<\sigma<1+\eta$, for any $\eta>0$, if $0<\alpha<1$, $\alpha\neq 1/2$, by works of Davenport and Heilbronn \cite{davenport1} and Cassels \cite{cassels}. Note that $L(1,1,s)=\zeta(s)$ and $L(1,1/2,s)=(2^s-1)\zeta(s)$, where $\zeta(s)$ is the Riemann zeta function.\\
When $0<\lambda <1$, we divide into three cases. If $\alpha$ is irrational, then Garunk\v{s}tis \cite{garunkstis2} has proved that $L(\lambda ,\alpha,s)$ has infinitely many zeros for $\sigma>1$ for every $\lambda $. If both $\alpha$ and $\lambda $ are rational, then $L(\lambda ,\alpha,s)$ has periodic coefficients, so it may be written as a linear combination over Dirichlet polynomials of Dirichlet $L$-functions associated with distinct primitive characters (see e.g. Saias and Weingartner \cite{saias}). It is easy to check that this linear combination is non trivial except for $\lambda =1/2$ and $\alpha=1/2$ or $\alpha=1$, for which one has $L(1/2,1,s)=(1-2^{1-s})\zeta(s)$ and $L(1/2,1/2,s)=2^sL(s,\chi)$, where $\chi$ is the non principal character mod $4$. Hence, excluding the two latter cases, by work of Saias and Weingartner \cite{saias} one gets that the Lerch zeta function $L(\lambda ,\alpha,s)$ has infinitely many zeros for $1<\sigma<1+\eta$, for any $\eta>0$. This result was already obtained by Laurin\v{c}ikas \cite{laurincikas2} in the particular case $\lambda =a/q$, $(a,q)=1$, $q$ odd prime, $\alpha=1$.  In this paper we deal with the remaining open case of irrational $\lambda $ and rational $\alpha$.

\begin{theorem}\label{theorem:lerch}
Let $0<\alpha,\lambda\leq 1$ be such that $\alpha$ is rational and $\lambda$ is irrational. Then, for any $\eta>0$, the Lerch zeta function $L(\lambda,\alpha,s)$ has infinitely many zeros for $1<\sigma<1+\eta$.
\end{theorem}

This result is actually an immediate consequence of the following more general one.

\begin{theorem}\label{theorem:zero_lin_twist}
Let $0<\lambda<1$ be an irrational number, and let $m$ and $k$ be integers such that $1\leq m\leq k$, $(m,k)=1$. Let $f(n)$ be a multiplicative arithmetic function and suppose that the associated Dirichlet series $F(s)=\sum_{n\geq 1} f(n)n^{-s}$ is absolutely convergent for $\sigma>1$. Suppose furthermore that
\begin{enumerate}[(i)]
\item\label{hp:PEP} for some positive integer $d$ we may write $\log F(s) = -\sum_p \sum_{j=1}^d\log(1-f_j(p)p^{-s})$ for $\sigma>1$, with $|f_j(p)|\leq 1$ for every prime $p$ and any $j=1,\ldots,d$;
\item\label{hp:RS} $\sum_{n\leq x} |f(n)|^2\ll x$ for every $x\geq 1$;
\item\label{hp:SW}  there exist positive constants $x_0$ and $A>6$ such that for any $x\geq x_0$ we have
$$\dsum_{p\leq x} |f(p)|^2 \gg \frac{x}{\log x}$$
and 
$$\dsum_{p\leq x} |f(p)|^2\chi(p) \ll  \frac{1}{\phi(q)}\frac{x}{\log^2 x}$$
for every Dirichlet character $\chi$ mod $q$, with $q\leq (\log x)^A$.
\end{enumerate}
Then, for any $\eta>0$, the twisted Dirichlet series $F(\lambda,m,k,s)=\sum_{n\equiv m\, (k)} f(n)\e{\lambda n}n^{-s}$ has infinitely many zeros for $1<\sigma<1+\eta$.
\end{theorem}

By the well known properties of the Dirichlet $L$-functions associated with primitive Dirichlet characters it is immediate to see that Theorem \ref{theorem:zero_lin_twist} may be applied to these $L$-functions and hence to $\zeta(s)$ (Theorem \ref{theorem:lerch}). We just note that in this case \ref{hp:SW} follows from Siegel--Walfisz's prime number theorem (see e.g. Davenport \cite[\S22]{davenport_mnt}).

Another example is given by $L$-functions associated with cusp forms of even weight for the full modular group since these functions are known to satisfy Ramanujan conjecture.
\begin{corollary*}
Let $0<\lambda<1$ be an irrational number, and let $m$ and $k$ be integers such that $1\leq m\leq k$, $(m,k)=1$. Let $g(z)=\sum_{n\geq 1} a(n)\e{nz}$ be a holomorphic cusp form of weight $\kappa\geq 2$ for the full modular group $\mathrm{SL}_2(\Z)$ which is an eigenfunction for all the Hecke operators with $a(1)=1$. Then, for any $\eta>0$, the Dirichlet series $\sum_{n\equiv m\, (k)} a(n)\e{\lambda n}n^{-s}$ has infinitely many zeros for $\frac{\kappa+1}{2}<\sigma<\frac{\kappa+1}{2}+\eta$.
\end{corollary*}
\begin{proof} It is well known that the multiplicativity of the function $a(n)$ comes from the fact that $g(z)$ is an eigenfunction for all the Hecke operators, as well as the polynomial Euler product of the associated $L$-function. Therefore, if we consider $f(n)=a(n)n^{-\frac{\kappa-1}{2}}$, then \ref{hp:PEP} follows with $d=2$ from Ramanujan conjecture, which was proved in this case by Deligne \cite[Theorem 8.2]{deligne1}. Furthermore, \ref{hp:RS} and \ref{hp:SW} are known to hold thanks to the properties of the Rankin--Selberg convolution: see e.g. Iwaniec and Kowalski \cite[\S14.9]{iwanieckowalski} for \ref{hp:RS}, and Perelli \cite{perelli0} and Ichihara \cite{ichihara} for \ref{hp:SW}. 
\end{proof}

Actually, conjecturally, the above properties should hold for all $L$-functions coming from irreducible unitary cuspidal automorphic representations on $\mathrm{GL}_r(\mathbb{A}_\Q)$, $r\geq 1$. Indeed \ref{hp:PEP} holds under Ramanujan's conjecture at every finite place; \ref{hp:RS} is known to hold by the properties of the Rankin--Selberg convolution; while \ref{hp:SW} would be a Siegel--Walfisz theorem for the Rankin--Selberg convolution, which may be obtained with standard means if it is the $L$-function of an automorphic form on $\mathrm{GL}_{2r}(\mathbb{A}_\Q)$, e.g. under Langlands' conjectures. Note that, for example, it is known that the Rankin--Selberg convolution of two cusp forms is the $L$-function of an automorphic form, usually cuspidal, on $\mathrm{GL}_4(\mathbb{A}_\Q)$ (see Ramakrishnan \cite{ramakrishnan}).

The proof of Theorem \ref{theorem:zero_lin_twist} follows the ideas of Saias and Weingartner \cite{saias} and will be presented in Section \ref{section:proof}. In Section \ref{section:exp_sums} we will prove a result on exponential sums with multiplicative coefficients over integers with small prime factors which generalizes Theorem 1 of Maier \cite{maier}; this is needed for the proof of Theorem \ref{theorem:zero_lin_twist}.
\medskip

We conclude with an application on the behavior of the the least upper bound of the real parts of the zeros of the linear twists of $L$-functions with respect to the twisting parameter $\lambda$. Namely, for any function $F(s)$ as in Theorem \ref{theorem:zero_lin_twist} and any $0<\lambda< 1$ let
$$\sigma^*(F,\lambda) = \inf\{\sigma_0\mid F(\lambda,1,1,s)\neq 0\hbox{ for }\sigma>\sigma_0\}.$$

\begin{theorem}\label{theorem:continuity} Let $F(s)$ be as in Theorem \ref{theorem:zero_lin_twist} and let $0<\lambda< 1$, $\lambda\neq 1/2$. Then for any sequence $\{x_n\}\subset(0,1]$ such that $x_n\rightarrow \lambda$ when $n\rightarrow \infty$ we have $\sigma^*(F,x_n)\rightarrow \sigma^*(F,\lambda)$ when $n\rightarrow \infty$.
\end{theorem}

\begin{remarks*}
(a) The proof of Theorem \ref{theorem:continuity} relies principally on the fact that $\sigma^*(F,\lambda)>1$, by Theorem \ref{theorem:zero_lin_twist} and by Theorem 3 of \cite{righetti}, and that $F(\lambda,1,1,\sigma)\not\rightarrow 0$ when $\sigma\rightarrow\infty$. In fact the proof can be adapted to work in much more generality: for example in our Ph.D. thesis \cite{righetti_phd} we showed that, for non-trivial linear combinations of Euler products, $\sigma^*$, as a function of the coefficients, is continuous outside of a precise Zariski closed set.\\
(b) By suitably modifying the proof of Theorem \ref{theorem:zero_lin_twist} one could obtain that $\sigma^*(F,a/q)-1\gg q^{-A}$ for some $A>0$ when $q\rightarrow\infty$.  However, this result shows that in reality $\sigma^*(F,a/q)$ depends only on $\norma{a/q}=\min(a/q,1-a/q)$, that is what one would expect.\\
(c) Concerning the distribution of the real parts of the zeros of $F(\lambda,m,k,s)$ for $\sigma>1$ we note that neither the proof of Theorem \ref{theorem:zero_lin_twist} nor the proof of Theorem \ref{theorem:continuity} yield that, when $\lambda$ is irrational, these real parts may be dense in some interval, although by Saias and Weingartner \cite{saias} and by \cite{righetti} we know that this is true when $\lambda$ is rational, except for the obvious cases.
\end{remarks*}

\subsection*{Notation}
For every Dirichlet character $\chi$ mod $q$, $q\geq 1$, we write
$$F_p(s,\chi) = \sum_{h\geq 0} \frac{f(p^h)\chi(p)^h}{p^{hs}}=\prod_{j=1}^d \left(1-\frac{f_j(p)\chi(p)}{p^s}\right)^{-1}\quad \hbox{and} \quad F(s,\chi)=\prod_p F_p(s,\chi)=\sum_{n\geq 1}\frac{f(n)\chi(n)}{n^{s}}.$$
For any Dirichlet series $L(s)=\sum_{n} a(n)n^{-s}$ and any completely multiplicative arithmetic function $\varphi(n)$, we write
$$L^\varphi(s)=\sum_{n}\frac{a(n)\varphi(n)}{n^{s}}.$$
We denote with $\phi(n)$ the Euler totient function and with $\tau(\chi)$ the Gauss' sum of any Dirichlet character $\chi$ mod $q$.
Moreover, we denote with $B_{R}(z)$ the closed disk in $\C$ of radius $R$ and center $z$.

%% file: exponential_sums.tex
\section{Exponential sums with multiplicative coefficients}\label{section:exp_sums}

In this section we collect some results on exponential sums with multiplicative coefficients both over the integers and over the integers with small prime factors (also known as \emph{smooth} or \emph{friable} numbers), which will be needed later in the proof of Theorem \ref{theorem:zero_lin_twist}.\\ 
Let $0<\alpha<1$ and let $\mathcal{F}_d$ be the class of all multiplicative arithmetic function $f(n)$ such that $|f(p)|\leq d$ for every prime $p$ and $\sum_{n\leq x} |f(n)|^2\leq d^2x$, $d>0$. First we consider the exponential sums
$$S(x,\alpha, f)=\sum_{n\leq x} f(n)\e{\alpha n}.$$
By work of Montgomery and Vaughan \cite{montvaughan}, we know that if $\alpha=a/q$, $(a,q)=1$, then we have (see Theorem 1 of \cite{montvaughan})
\begin{equation}\label{eq:S_rational}
S(x,a/q, f)\ll_d \frac{x}{\log 2x}+\frac{x}{\sqrt{\phi(q)}}+\sqrt{qx}(\log(2x/q))^{3/2},\quad \hbox{if }x\geq q,\quad \hbox{uniformly for }f\in\mathcal{F}_d.
\end{equation}
If $\alpha$ is irrational and $a/q$ is such that $(a,q)=1$ and $|\alpha -a/q|\leq q^{-2}$, we have (see \cite[\S6]{montvaughan})
\begin{equation}\label{eq:S_irrational_large_r}
S(x,\alpha, f)\ll_d \frac{x}{\log x}+\frac{x}{\sqrt{q}}(\log q)^{3/2},\quad \hbox{if }x\geq q^2,\quad \hbox{uniformly for }f\in\mathcal{F}_d.
\end{equation}

%
Now, let $P(n)$ be the largest prime factor of an integer $n$, with the convention $P(1)=1$. Then we consider the exponential sums
$$S(x,y,\alpha, f)=\sum_{\substack{n\leq x\\P(n)\leq y}} f(n)\e{\alpha n},\quad 1\leq y\leq x.$$
When $\alpha=s/r$, with $(s,r)=1$ and $r$ prime, and $f(n)$ is multiplicative with $|f(n)|\leq 1$ for every $n$, there is a result of Maier \cite[Theorem 1]{maier}: for every $\eps_0>0$ and every $A>0$ we have
\begin{equation}\label{eq:maier}
S(x,y,s/r, f)\ll_{\eps_0,A} \frac{\Psi(x,y)}{\sqrt{r}},
\end{equation}
where $\Psi(x,y)=|\{n\leq x: P(n)\leq y\}|$, uniformly for $\exp{(\log x)^{\eps_0}}< y\leq x$, $r\leq (\log x)^A$ and $f(n)$ multiplicative such that $|f(n)|\leq 1$. We need an analogue of \eqref{eq:maier} for $f\in\mathcal{F}_d$. To this purpose we prove the following general result.
\begin{theorem}\label{theorem:ub_smooth}
Let $f(n)$ be a multiplicative function. Let $\eps_0>0$ be arbitrary small, $A>0$ be arbitrary large, $\exp{(\log x)^{\eps_0}}<y\leq x$, $r\leq (\log x)^A$ be a prime, $(s,r)=1$. Then we have
\begin{equation}\label{eq:mamova}
S(x,y,s/r, f)\ll_{\eps_0,A} \left(\sum_{\substack{n\leq x\\P(n)\leq y}} |f(n)|^2\right)^\frac{1}{2}\left(\frac{\Psi(x,y)}{\sqrt{r}}\right)^\frac12 + |f(r)|\left(\sum_{\substack{n\leq x/r\\P(n)\leq y}} |f(n)|^2\right)^{\frac12}\Psi(x/r,y)^{\frac12},
\end{equation}
for any $f(n)$ as above.
\end{theorem}
\begin{remarks*}(a) This result is not optimal, but it suffices for the estimates needed in the proof of Theorem \ref{theorem:zero_lin_twist}.\\
(b) The second term in \eqref{eq:mamova} may be removed if $y$ (or $r$) is sufficiently large, as one may see by combining Theorems 2 and 3 of Hildebrand and Tenenbaum \cite{hildten86}. \\
(c) It is possible to similarly adapt the arguments of de la Bret\`eche and Tenenbaum \cite[Theorem 2.1]{dlbten2} to obtain analogous results with $r$ not necessarily prime or $\alpha$ irrational.
\end{remarks*}
\begin{proof} The proof is a straightforward adaptation of Theorem 1 of Maier \cite{maier} by suitably using Cauchy-Schwarz inequality. Therefore here we sketch only the necessary details and refer to \cite{maier} for the rest. For easier reference we keep the same notation as in \cite{maier}. We hence write
\begin{equation}\label{eq:3_pezzi}
\sum_{\substack{n\leq x\\P(n)\leq y}} f(n)\e{\frac{sn}{r}} = \sum_{n\in\mathfrak{m}_1}f(n)\e{\frac{sn}{r}} + \sum_{n\in\mathfrak{m}_2}f(n)\e{\frac{sn}{r}} + \sum_{\substack{n\leq x\\P(n)\leq y\\ r|n}} f(n)\e{\frac{sn}{r}},
\end{equation}
where $\mathfrak{m}_1$ and $\mathfrak{m}_2$ are defined in \cite[p. 212]{maier}.\\
For the last term in \eqref{eq:3_pezzi} we use Cauchy-Schwarz inequality and we obtain
$$\abs{\sum_{\substack{n\leq x\\P(n)\leq y\\ r|n}} f(n)\e{\frac{sn}{r}}}\leq d\left(\sum_{\substack{m\leq x/r\\P(m)\leq y}} |f(m)|^2\right)^{\frac12}\Psi(x/r,y)^{\frac12}.$$
For the second term in \eqref{eq:3_pezzi} we recall that from the proof of Lemma 3.9 \cite[p. 217]{maier} we have $|\mathfrak{m}_2|\ll \Psi(x,y)/M_0$, where $M_0=\exp{(\log x)^{\eps_0/2}}\gg (\log x)^{A/2}\gg \sqrt{r}$. We then use Cauchy-Schwarz inequality on the second term and we obtain
$$\abs{ \sum_{n\in\mathfrak{m}_2}f(n)\e{\frac{sn}{r}} }\leq \left(\sum_{\substack{n\leq x\\P(n)\leq y}} |f(n)|^2\right)^{\frac12} |\mathfrak{m}_2|^{\frac12}\ll \left(\sum_{\substack{n\leq x\\P(n)\leq y}} |f(n)|^2\right)^\frac{1}{2}\left(\frac{\Psi(x,y)}{\sqrt{r}}\right)^\frac12.$$
We split the first term in \eqref{eq:3_pezzi} in two parts, namely
\begin{equation}\label{eq:2_pezzi}
\sum_{n\in\mathfrak{m}_1}f(n)\e{\frac{sn}{r}} = \sum_{n\in\mathfrak{m}_1^{(*)}}f(n)\e{\frac{sn}{r}} + \sum_{n\in\mathfrak{m}_1-\mathfrak{m}_1^{(*)}}f(n)\e{\frac{sn}{r}},
\end{equation}
where $\mathfrak{m}_1^{(*)}$ is defined in \cite[Definition 2.5]{maier}. For the first term of \eqref{eq:2_pezzi} we look at the proof of \cite[Lemma 3.8]{maier}. Due to the definition of $\mathfrak{m}_1$, $\mathfrak{m}_1^{(1)}$ and $\mathfrak{m}_1^{(2)}$ we may write
$$\sum_{n\in\mathfrak{m}_1^{(*)}}f(n)\e{\frac{sn}{r}} = \sum_{m_0\leq M_0} \sum_{\substack{\overline{l}\, such\, that\\\mathfrak{m}_{1,\overline{l},m_0}\,proper}}\sum_{n\in\mathfrak{m}_{1,\overline{l},m_0}} f(n)\e{\frac{sn}{r}},$$
where $\mathfrak{m}_{1,\overline{l},m_0}$ is defined in \cite[p. 212]{maier}, as well as the property of being \emph{proper}.  The inner sum is called $Q_2$ in \cite[p. 216]{maier}, and may be written as a double sum by decomposing a proper $\mathfrak{m}_{1,\overline{l},m_0}$ as the product set (given by all the products) $m_0\cdot S_1\cdot S_2$, where $S_1$ and $S_2$ are defined in \cite[p. 215]{maier} through Definition 2.6 of \cite{maier}. By using three times Cauchy-Schwarz inequality and the multiplicativity of $f(n)$ we obtain (cf. \cite[p. 216]{maier})
\begin{spliteq*}
|Q_2| \leq& |f(m_0)|\left(\sum_{n_1\in S_1}|f(n_1)|^2\right)^\frac12\left(\sum_{n_2^{(1)}\in S_2} f(n_2^{(1)})\sum_{n_2^{(2)}\in S_2}\conj{f(n_2^{(2)})}\sum_{n_1\in S_1}\e{(n_2^{(1)}-n_2^{(2)})\frac{m_0s}{r}n_1}\right)^\frac12\\
\leq & |f(m_0)|\left(\sum_{n_1\in S_1}|f(n_1)|^2\right)^\frac12\left(\left(\sum_{n_2\in S_2}|f(n_2)|^2\right)^\frac12 \cdot\sum_{n_2^{(1)}\in S_2} |f(n_2^{(1)})|\left(\sum_{n_2^{(2)}\in S_2}|Q_3|^2\right)^\frac12\right)^\frac12\\
\leq & \left(\sum_{n\in\mathfrak{m}_{1,\overline{l},m_0}} |f(n)|^2\right)^\frac12\left(\sum_{n_2^{(1)}\in S_2} \sum_{n_2^{(2)}\in S_2}|Q_3|^2\right)^\frac14.
\end{spliteq*}
With the same arguments as in \cite[p. 216]{maier},\footnote{In \cite[p. 216]{maier} there is a misprint in the last inequality: the last $S_1$ should be $S_2$.} we obtain
$$|Q_2|\ll \left(\sum_{n\in\mathfrak{m}_{1,\overline{l},m_0}} |f(n)|^2\right)^\frac12\left(|S_1|^2\frac{|S_2|^2}{r} + \frac{|S_1|^2}{(r-1)^{2}} |S_2|^2\right)^\frac14\ll \left(\sum_{n\in\mathfrak{m}_{1,\overline{l},m_0}} |f(n)|^2\right)^\frac12 |\mathfrak{m}_{1,\overline{l},m_0}|^\frac12\, r^{-\frac14}.$$
Therefore, using again twice Cauchy-Schwarz inequality, we get
\begin{spliteq*}
\sum_{n\in\mathfrak{m}_1^{(*)}}f(n)\e{\frac{sn}{r}} \ll &\, r^{-\frac14}\sum_{m_0\leq M_0} \left(\sum_{\substack{\overline{l}\, such\, that\\\mathfrak{m}_{1,\overline{l},m_0}\,proper}}\sum_{n\in\mathfrak{m}_{1,\overline{l},m_0}} |f(n)|^2\right)^\frac12\left(\sum_{\substack{\overline{l}\, such\, that\\\mathfrak{m}_{1,\overline{l},m_0}\,proper}}|\mathfrak{m}_{1,\overline{l},m_0}|\right)^\frac12\\
\ll &\, r^{-\frac14}\left(\sum_{\substack{n\leq x\\P(n)\leq y}} |f(n)|^2\right)^\frac12\Psi(x,y)^\frac12.
\end{spliteq*}
For the second term in \eqref{eq:2_pezzi} we use Cauchy-Schwarz inequality and, by Lemma 3.5 of \cite{maier}, we obtain
$$\abs{\sum_{n\in\mathfrak{m}_1-\mathfrak{m}_1^{(*)}}f(n)\e{\frac{sn}{r}}}\leq \left(\sum_{\substack{n\leq x\\P(n)\leq y}} |f(n)|^2\right)^{\frac12} \abs{\mathfrak{m}_1-\mathfrak{m}_1^{(*)}}^{\frac12}\ll \left(\sum_{\substack{n\leq x\\P(n)\leq y}} |f(n)|^2\right)^{\frac12} \left(\frac{\Psi(x,y)}{(\log y)^A}\right)^\frac12.$$
The result then follows from the fact that $(\log x)^A\geq \sqrt{r}$ and by repeating the proof with $A=A/\eps_0$.
\end{proof}

%% file: existence_zeros.tex
\section{Proof of Theorem \ref{theorem:zero_lin_twist}}\label{section:proof}

We fix $0<\delta<1/3$ such that $(1+\delta)A>8$ and, by work of many authors, from Vinogradov \cite{vinogradov} to Matom\"aki \cite{matomaki}, we know that there is an infinite sequence of rational numbers $a_n/q_n$ such that $q_n$ is prime, $(a_n,q_n)=1$ and 
\begin{equation}\label{eq:dioph}
\abs{\lambda-\frac{a_n}{q_n}}<\frac{1}{q_n^{1+\delta}}.
\end{equation}
We hence fix an integer $\ell \geq 2$ such that $(\ell,k)=1$ and $f(\ell)\neq 0$, we take a fixed arbitrarily large prime number $q>k+\ell$ such that \eqref{eq:dioph} holds for some coprime number $a<q$, and we set $Q=q^{(1+\delta)/8}$. Then, if $\psi$ runs among the Dirichlet characters mod $k$ and $\chi$ runs among those mod $q$, since $(m,k)=1$ and $(a,q)=1$, by the orthogonality of characters we get 
\begin{spliteq*}
F(\lambda,m,k,s)=& \sum_{n\equiv m\,(k)} \frac{f(n)\e{an/q}}{n^s}+\sum_{n\equiv m\,(k)} \frac{f(n)[\e{\lambda  n}-\e{an/q}]}{n^s}\\
=& \frac{1}{\phi(k)}\sum_{\psi} \conj{\psi(m)}\sum_{b=0}^{q-1} \e{\frac{ab}{q}}\sum_{n\equiv b\,(q)} \frac{f(n)\psi(n)}{n^s}+R(s,q)\\
=& \frac{1}{\phi(k)\phi(q)}\sum_{\psi} \conj{\psi(m)} \sum_{\chi}\left(\sum_{b=1}^{q-1} \e{\frac{ab}{q}}\conj{\chi}(b)\right)F(s,\psi\chi)\\
&\quad+\frac{1}{\phi(k)}\sum_{\psi} \conj{\psi(m)}[F(s,\psi)-F(s,\psi\chi_0)]+R(s,q)\\
=& \frac{1}{\phi(kq)}\sum_{\psi} \sum_{\chi\neq \chi_0}\conj{\psi(m)}\chi(a)\tau(\conj{\chi})F(s,\psi\chi)\\
&\quad-\frac{1}{\phi(kq)}\sum_{\psi} \conj{\psi(m)}\left[1-\phi(q)(F_q(s,\psi)-1)\right]F(s,\psi\chi_0)+R(s,q),
\end{spliteq*}
where $\chi_0$ is the principal character mod $q$ and
$$R(s,q)=\sum_{n\equiv m\,(k)} \frac{f(n)[\e{\lambda  n}-\e{an/q}]}{n^s}.$$
By Bohr's equivalence theorem, solving the equation $F(\lambda ,m,k,\sigma+it)=0$ is equivalent to finding a completely multiplicative function $\varphi(n)$, with $|\varphi(p)|=1$ for every prime $p$, such that $F^\varphi(\lambda ,m,k,\sigma)=0$ (cf. e.g. Chapter 8 of Apostol \cite{apostol} or the introduction of \cite{righetti}). For convenience we write $\varphi(p)=p^{-it_p}$, with $t_p\in\R$ to be determined for every prime $p$.\\
We denote with $\chi_1,\ldots,\chi_{q-2}$ the $q-2$ primitive Dirichlet characters mod $q$ in such a way that $\chi_h(n)=\chi_{n-1}(h+1)$ for every $h=1,\ldots,q-2$, $n=2,\ldots,q-1$; this may be achieved by taking, for example, $\chi_h(n)=\e{\frac{\nu(h+1)\nu(n)}{q-1}}$, where $\nu(n)$ is the index of $n$ relative to a fixed primitive root mod $q$ (cf. e.g. Davenport \cite[p. 29]{davenport_mnt}). Then we consider the $\phi(kq)$ functions
$$X_{\psi,j}(s,q)=\left(1-W_{\psi,j}(s,q)\right)\prod_{p\leq \exp{Q}}F_p(\sigma,\psi\chi_j),\quad \psi\hbox{ mod }k,\, j=0,\ldots,q-2,$$
where $W_{\psi,j}(s,q)=0$ for $j=0$ or $j> Q^2$ and any $\psi\hbox{ mod }k$, and
\begin{equation}\label{eq:W_psi_j}
W_{\psi,j}(s,q) = \frac{\conj{\psi(\ell)}\conj{\chi_j(\ell)}}{\conj{\psi(m)}\chi_j(a)\tau(\conj{\chi_j})}\frac{R(s,q)+ \dfrac{1}{\phi(kq)}\dsum_{\psi} \conj{\psi(m)}\sum_{h=0}^{q-2}\chi_h(a)\tau(\conj{\chi_h})\prod_{p\leq \exp{Q}}F_p(\sigma,\psi\chi_h)}{\dfrac{1}{\phi(kq)}\dsum_{\psi} \conj{\psi(\ell)}\dsum_{h=1}^{Q^2}\conj{\chi_h(\ell)}\prod_{p\leq \exp{Q}}F_p(\sigma,\psi\chi_h)},
\end{equation}
for $1\leq j\leq Q^2$ and any $\psi\hbox{ mod }k$. It is easy to check that 
$$ \frac{1}{\phi(kq)}\sum_{\psi}\sum_{j=1}^{q-2}\conj{\psi(m)}\chi_j(a)\tau(\conj{\chi_j})X_{\psi,j}(s,q)-\frac{1}{\phi(kq)}\sum_{\psi} \conj{\psi(m)}X_{\psi,0}(s,q)+R(s,q)=0$$
identically for $\sigma>1$. Hence, the result would follow if we can find $\sigma>1$ such that the system of equations
\begin{equation}\label{eq:system}
\left\{\begin{array}{ll}
\left[1-(q-1)(F_q(\sigma+it_q,\psi)-1)\right]\dprod_{p> \exp{Q}}F_p(\sigma+it_p,\psi)=1& \psi\hbox{ mod }k,\\
\dprod_{p> \exp{Q}}F_p(\sigma+it_p,\psi\chi_j)=1-W_{\psi,j}^\varphi(\sigma,q), & \psi\hbox{ mod }k,\, j=1,\ldots,q-2,
\end{array}\right.
\end{equation}
has a solution $\{t_p\}_{p}$, with $t_p=0$ for every prime $p\leq \exp{Q}$, and $p\neq q$ if $f(q)\neq 0$, in which case $t_q=(\arg(f(q))-2\pi /q - \pi)/\log q$. Note that $F_q(s,\psi)^{-1}$ is a polynomial of degree $d$ in the variable $q^{-s}$, and hence $1-(q-1)(F_q(s,\psi)-1)$  vanishes only on at most $d$ vertical lines. So we may suppose, without loss of generality, that $\sigma$ is not any of these $d\cdot \phi(k)$ values. Therefore, we may take the principal branch of the logarithm and reduce \eqref{eq:system} to the equivalent system
\begin{equation*}
\left\{\begin{array}{ll}
\dsum_{p>\exp{Q}}\log F_p(\sigma+it_p,\psi)=-\log(1-(q-1)(F_q(\sigma+it_q,\psi)-1)) & \psi\hbox{ mod }k,\\
\dsum_{p>\exp{Q}}\log F_p(\sigma+it_p,\psi\chi_j)=\log (1-W_{\psi,j}^\varphi(\sigma,q)), & \psi\hbox{ mod }k,\, j=1,\ldots,q-2.
\end{array}\right.
\end{equation*}
Actually, it is better to work with another equivalent system, obtained by multiplication with the inverse of the matrix
$$(\psi(b)\chi_j(b))_{\substack{\psi\,\mathrm{ mod }\,k, \,j=0,\ldots,q-2\\b=1,\ldots,kq,\,(b,kq)=1}},$$
i.e. with the system
\begin{equation}\label{eq:log_system}
\sum_{j=1}^d\sum_{h=1}^{\infty}\sum_{\substack{p> \exp{Q}\\p^h\equiv b\,(kq)}}\frac{f_j(p)^h}{hp^{h(\sigma+it_p)}}=Y_b^\varphi(\sigma,q),\quad b=1,\ldots,kq,\, (b,kq)=1,
\end{equation}
where
$$Y_b^\varphi(\sigma,q)=\frac{1}{\phi(kq)}\sum_{\psi}\sum_{j=1}^{Q^2}\conj{\psi(b)}\conj{\chi_j(b)}\log (1-W_{\psi,j}^\varphi(\sigma,q))-\frac{1}{\phi(kq)}\sum_{\psi}\conj{\psi(b)}\log\!\left(1-(q-1)\sum_{h=1}^\infty \frac{f(q^h)\psi(q^h)}{q^{h(\sigma+it_q)}}\right)\!,$$
$b=1,\ldots,kq$, $(b,kq)=1$.

Since it is sufficient to solve this system when $q$ is  sufficiently large, we study the behavior of the LHS and of the RHS of \eqref{eq:log_system} as $q\rightarrow \infty$ ($a$ has to be considered as a function of $q$ so that \eqref{eq:dioph} holds). For the LHS we note that from \ref{hp:PEP}, \ref{hp:SW}, the orthogonality of characters and the classical prime number theorem (cf. e.g. \cite[\S18]{davenport_mnt}), we obtain that there exist positive constants $B$, $C$ and $D$, such that for any $b$ and $q$, with $(b,kq)=1$, $q\geq B$, and for any $x\geq \exp{q^{1/A}}$ we have
\begin{equation*}\label{eq:SW}
\frac{C}{\phi(q)}\frac{x}{\log x}\leq \dsum_{\substack{p\leq x\\p\equiv b\, (kq)}} |f(p)|^2 \leq  \frac{D}{\phi(q)}\frac{x}{\log x}.
\end{equation*}
Then, for the choice of $\delta$ and $Q$ we obtain 
\begin{equation}\label{eq:lb_main_term}
\sum_{\substack{p>\exp{Q}\\p\equiv b\,(kq)}}\frac{|f(p)|}{p^{\sigma}} \geq\frac{1}{d}\sum_{\substack{p> \exp{Q}\\p\equiv b\,(kq)}}\frac{|f(p)|^2}{p^{\sigma}} \geq  -\frac{D\exp{(1-\sigma)Q}}{d\phi(kq)Q}+\frac{C\sigma}{d\phi(kq)}\int_{ \exp{Q}}^\infty \frac{dx}{x^\sigma\log x}\geq -\frac{C_1}{qQ}+\frac{C_2}{q}\int_{ \exp{Q}}^\infty \frac{dx}{x^\sigma\log x},
\end{equation}
for some positive constants $C_1$ and $C_2$, for every $1<\sigma\leq 1+\eta$, if $q$ is sufficiently large, $b=1,\ldots,kq$, $(b,kq)=1$. Moreover, by \ref{hp:PEP},  
Chebyshev's bounds for the prime-counting function (see e.g. Davenport \cite[\S7(1)]{davenport_mnt}) and a well known bound for the exponential integral function (see e.g. Abramowitz and Stegun \cite[5.1.20]{abramowitzstegun}), we get
\begin{spliteq}\label{eq:ub_main_error}
\abs{\sum_{j=1}^d\sum_{h=2}^{\infty}\sum_{\substack{p> \exp{Q}\\p^h\equiv b\,(kq)}}\frac{f_j(p)^h}{hp^{h(\sigma+it_p)}} }& \leq d\sum_{h=2}^{\infty}\sum_{p> \exp{Q}} \frac{1}{hp^{h\sigma}} \ll\sum_{h=2}^{\infty}\int_{\exp{Q}}^\infty \frac{x^{2-h}}{x^2\log x}dx\\
&\ll \int_{Q}^\infty \frac{\exp{-w}}{w}dw\cdot\sum_{h=2}^{\infty}\exp{(2-h)Q}\ll \frac{\exp{-Q}}{Q}\ll\frac{1}{q},
\end{spliteq}
uniformly for $1\leq \sigma\leq 1+\eta$, if $q$ is sufficiently large, $b=1,\ldots,kq$, $(b,kq)=1$.

For the RHS we have
\begin{equation}\label{eq:ub_q}
\log\!\left(1-(q-1)\sum_{h=1}^\infty \frac{f(q^h)\psi(q)^h}{q^{h(\sigma+it_q)}}\right)=\log\!\left(1+q^{1-\sigma}|f(q)|\psi(q)\exp{2\pi i/q}+\Og{1/q}\right)\ll 1,
\end{equation}
uniformly in $q$ for every $\sigma\geq 1$ and every $\psi$ mod $k$, since $\psi(q)\exp{2\pi i /q}\neq -1$ and hence the argument is bounded away from 0 independently from $q$. Moreover, we have 
\begin{equation}\label{eq:first_term_log}
\log (1-W_{\psi,j}^\varphi(\sigma,q))=-W_{\psi,j}^\varphi(\sigma,q)+O\!\left(\frac{\log^2 q}{q}\right),\quad \psi\hbox{ mod }k, \, 1\leq j\leq Q^2,
\end{equation}
uniformly for $Q^{-3/2} \leq \sigma-1\leq Q^{-1}$, if $q$ is sufficiently large. Indeed, by \ref{hp:RS}, \eqref{eq:dioph}, Cauchy-Schwarz inequality and summation by parts we have that
\begin{equation}\label{eq:small_terms}
\sum_{\substack{ n\equiv m\,(k)\\n\leq q^{1+\delta}}}\frac{|f(n)||\e{an/q}-\e{\lambda n}|}{n^\sigma}\ll \log q,
\end{equation}
uniformly for $\sigma\geq 1$. Moreover, let $l/r$ be a convergent to the continued fraction of $\lambda$, then by Matom\"aki \cite{matomaki} there exists at least one prime $q$ such that \eqref{eq:dioph} holds for some $(a,q)=1$ and such that $r^{2/(1+\delta)}< q \leq 2 r^{2/(1+\delta)}$, if $r$ is sufficiently large (depending on $\delta$). Therefore, since $|\lambda-l/r|<r^{-2}$, \eqref{eq:S_irrational_large_r} becomes
\begin{equation}\label{eq:S_irrational_large}
S(x,\lambda ,f)\ll \frac{x}{\log x}+\frac{x}{q^{(1+\delta)/4}}(\log q)^{3/2},\quad \hbox{if }x\geq q^{1+\delta},\quad \hbox{uniformly for }f\in\mathcal{F}_d.
\end{equation}
Then, by \ref{hp:RS}, \eqref{eq:S_rational} and \eqref{eq:S_irrational_large} we get
\begin{spliteq}\label{eq:bound_int_S}
\int_{q^{1+\delta}}^\infty &\frac{|S(x,\lambda,f\cdot\psi\cdot\varphi)|+|S(x,a/q,f\cdot\psi\cdot\varphi)|}{x^{\sigma+1}}dx\\
 &\qquad\qquad\ll \int_{q^{1+\delta}}^\infty\left[\frac{1}{x^{\sigma}\log x}+\frac{(\log q)^{3/2}}{Q^2x^{\sigma}}+\frac{\sqrt{q}(\log(2x/ q))^{3/2}}{x^{\sigma+1/2}}\right]dx \\
&\qquad\qquad \ll \int_{q^{1+\delta}}^{ \exp{Q^2}}\frac{dx}{x\log x}+\int_{\exp{Q^2}}^{\infty}\frac{dx}{x^{\sigma}\log x}+\frac{(\log q)^{3/2}}{Q^2}\frac{q^{(1+\delta)(1-\sigma)}}{\sigma-1}+\int_{ q^\delta}^\infty \frac{\log(2x)^{3/2}}{x^{3/2}}dx\\
&\qquad\qquad\ll  \log q,
\end{spliteq}
uniformly for $Q^{-3/2}\leq \sigma-1\leq \eta$, if $r$, and thus $q$, is sufficiently large. Hence, summing by parts, by \eqref{eq:dioph}, \eqref{eq:small_terms} and \eqref{eq:bound_int_S}, we obtain
\begin{equation}\label{eq:bound_R}
|R^\varphi(\sigma,q)| \ll \log q
\end{equation}
uniformly for $Q^{-3/2} \leq \sigma-1\leq \eta$, if $r$, and thus $q$, is sufficiently large. 

Moreover, by suitably adapting the proof of Theorem 11 of Fouvry and Tenenbaum \cite{fouvrytenenbaum}, we obtain that
\begin{equation}\label{eq:bound_smooth}
\frac{1}{\phi(kq)}\dsum_{\psi} \conj{\psi(m)}\sum_{h=0}^{q-2}\chi_h(a)\tau(\conj{\chi_h})\prod_{p\leq \exp{Q}}F_p(\sigma,\psi\chi_h)\ll \log q,
\end{equation}
uniformly for $1\leq \sigma\leq 1+\eta$, if $q$ is sufficiently large. To prove this we recall the well known bound on $\Psi(x,y)$ due to de Bruijn \cite[(1.9)]{debruijn}, i.e.
\begin{equation}\label{eq:debruijn}\Psi(x, \exp{Q})\ll x^{1-c/Q},\quad \hbox{for a positive constant $c$, uniformly for $x>1$ and $q\geq 1$.}
\end{equation}
Furthermore, we note that, by the absolute convergence for $\sigma>1$, it is enough to prove \eqref{eq:bound_smooth} for $\sigma=1$. Since for any $\psi$ mod $k$ we have
\begin{spliteq*}
\frac{1}{\phi(q)}\sum_{h=0}^{q-2}\chi_h(a)\tau(\conj{\chi_h})\prod_{p\leq \exp{Q}}F_p(\sigma,\psi\chi_h)& = \sum_{\substack{n\geq 1\\ P(n)\leq \exp{Q}}}\frac{f(n)\varphi(n)\psi(n)\chi_0(n)\e{an/q}}{n^\sigma}\\
&= \sigma\int_1^\infty S(x,\exp{Q},a/q,f\cdot \varphi\cdot\psi\cdot \chi_0)\frac{dx}{x^{\sigma+1}},
\end{spliteq*}
to prove \eqref{eq:bound_smooth} we split the latter integral into four parts and bound them separately. For starters we trivially have
$$\int_1^q |S(x,\exp{Q},a/q,f\cdot \psi\cdot \chi_0)|\frac{dx}{x^{\sigma+1}}\leq \log q,$$
uniformly for $\sigma\geq 1$. Since $f\cdot\varphi\cdot \psi\cdot\chi_0\in\mathcal{F}_d$ by \ref{hp:PEP} and \ref{hp:RS}, by \eqref{eq:S_rational} we get
$$\int_q^{\exp{Q}} |S(x,\exp{Q},a/q,f\cdot \psi\cdot \chi_0)|\frac{dx}{x^{\sigma+1}}\ll \int_q^{\exp{Q}} \left[\frac{1}{x\log x}+\frac{1}{x\sqrt{q}}+\frac{(\log(x/q))^{3/2}}{q(x/q)^{3/2}}\right]dx\ll \log q + \frac{Q}{\sqrt{q}} + 1,$$
uniformly for $\sigma\geq 1$. Moreover, by Theorem \ref{theorem:ub_smooth} with $A=8/(1+\delta)$ and $\eps_0=2/3$, since $q$ is prime, we get
$$\int_{\exp{Q}}^{\exp{Q^{3/2}}} |S(x,\exp{Q},a/q,f\cdot \psi\cdot \chi_0)|\frac{dx}{x^{\sigma+1}}\ll \int_{\exp{Q}}^{\exp{Q^{3/2}}} \frac{dx}{xq^{1/4}}\ll \frac{Q^{3/2}}{q^{1/4}},$$
uniformly for $\sigma\geq 1$. Finally, by \eqref{eq:debruijn}, we have
$$\int_{\exp{Q^{3/2}}}^{\infty} |S(x,\exp{Q},a/q,f\cdot \psi\cdot \chi_0)|\frac{dx}{x^{\sigma+1}}\leq \int_{\exp{Q^{3/2}}}^{\infty} x^{-1-c/Q}dx \leq \exp{-c\sqrt{Q}}Q^{3/2},$$
uniformly for $\sigma\geq 1$. Putting these together and summing over $\psi$ mod $k$ we obtain that \eqref{eq:bound_smooth} holds uniformly for $1\leq \sigma\leq 1+\eta$, if $q$ is sufficiently large.

Finally, we note that by the orthogonality of Dirichlet characters we have
\begin{spliteq*}
\dfrac{1}{\phi(kq)}\dsum_{\psi} &\conj{\psi(\ell)}\dsum_{h=1}^{Q^2}\conj{\chi_h(\ell)}\prod_{p\leq \exp{Q}}F_p(\sigma,\psi\chi_h)\\
&=\sum_{\substack{n\equiv \ell\,(kq)\\ P(n)\leq \exp{Q}}}\frac{f(n)}{n^\sigma}-\frac{1}{q-1}\sum_{\substack{n\equiv \ell\,(k)\\ P(n)\leq \exp{Q}}}\frac{f(n)}{n^\sigma}\left[\chi_0(n)+\sum_{Q^2<h\leq q-2}\chi_h(n)\conj{\chi_h(\ell)}\right].
\end{spliteq*}
By Cauchy-Schwarz inequality and \ref{hp:RS}, we have
$$\sum_{\substack{n\leq x\\ n\equiv \ell\, (kq)}}|f(n)|\leq \left(\sum_{n\leq x}|f(n)|^2\right)^\frac12 \left(\sum_{\substack{n\leq x\\ n\equiv \ell\, (kq)}}1\right)^\frac12\ll \frac{x}{\sqrt{q}},\quad \hbox{for every }x, \,q> k+\ell. $$
Hence, by the triangle inequality, partial summation and the P\'olya-Vinogradov inequality (cf. e.g. Davenport \cite[\S23]{davenport_mnt}), we obtain
\begin{spliteq*}
\Bigg|\sum_{\substack{n\equiv \ell\, (kq)\\ P(n)\leq\exp{Q}\\ n\geq kq}}\frac{f(n)}{n^{\sigma}}-\frac{1}{q-1}\sum_{\substack{n\equiv \ell\,(k)\\ P(n)\leq \exp{Q}}} &\frac{f(n)}{n^\sigma}\Bigg[\chi_0(n)+\sum_{Q^2<h\leq q-2}\chi_h(n)\conj{\chi_h(\ell)}\Bigg]\Bigg| \\
&\leq \sum_{\substack{n\equiv \ell\, (kq)\\n\geq kq}}\frac{|f(n)|}{n^{\sigma}}+\frac{1}{q-1}\sum_{\substack{n\equiv \ell\,(k)\\(n,q)=1}} \frac{|f(n)|}{n^{\sigma}}\abs{1+\sum_{Q^2+1<h\leq q-1}\chi_{n-1}(h)\conj{\chi_{\ell-1}(h)}} \\
&\ll \frac{\sigma}{\sqrt{q}}\int_{kq}^\infty \frac{dx}{x^\sigma}+ \frac{\sigma\log q}{\sqrt{q}}\int_1^\infty \frac{dx}{x^\sigma}\ll \frac{Q^{3/2}\log q}{\sqrt{q}},
\end{spliteq*}
uniformly for $Q^{-3/2}\leq \sigma-1 \leq \eta$, if $q$ is sufficiently large. Therefore by the triangle inequality we get
\begin{equation}\label{eq:lb_denominator}
\abs{\dfrac{1}{\phi(kq)}\dsum_{\psi} \conj{\psi(\ell)}\dsum_{h=1}^{Q^2}\conj{\chi_h(\ell)}\prod_{p\leq \exp{Q}}F_p(\sigma,\psi\chi_h)}\geq \frac{1}{2}\frac{|f(\ell)|}{\ell^{1+\eta}},
\end{equation}
uniformly for $Q^{-3/2}\leq \sigma-1 \leq \eta$, if $q$ is sufficiently large. Putting together \eqref{eq:bound_R}, \eqref{eq:bound_smooth}, \eqref{eq:lb_denominator} and the well known fact that $|\tau(\chi)|=\sqrt{q}$ for any primitive Dirichlet character $\chi$ mod $q$ (see e.g. Davenport \cite[\S9(5)]{davenport_mnt}), we obtain
$$W_{\psi,j}^\varphi(\sigma,q)\ll\frac{\log q}{\sqrt{q}}, \qquad \psi\hbox{ mod }k, \, 1\leq j\leq Q^2,$$
uniformly for $Q^{-3/2}\leq \sigma-1\leq \eta$, if $r$, and thus $q$, is sufficiently large. This proves \eqref{eq:first_term_log}.\\
It is then easy to check that, by \eqref{eq:ub_q} and \eqref{eq:first_term_log} we trivially have
\begin{equation}\label{eq:ub_Y_b}
Y_b^\varphi(\sigma,q)\ll \frac{1}{q}\left(1+\frac{Q^2}{\sqrt{q}}\log q\right)\ll \frac{1}{q},\qquad b=1,\ldots,kq,\, (b,kq)=1,
\end{equation}
uniformly for $Q^{-3/2} \leq \sigma-1\leq \eta$, if $r$, and thus $q$, is sufficiently large.

Now, we want to use Saias and Weingartner's approach to find solutions of \eqref{eq:log_system} with Brouwer fixed point theorem (see \cite{saias}). We hence rewrite \eqref{eq:log_system} as
\begin{equation}\label{eq:sw_system}
\sum_{\substack{p> \exp{Q}\\p\equiv b\,(kq)}}\frac{f(p)}{p^{\sigma+it_p}}-E_b(\sigma,(t_p)_{p>\exp{Q}})=0,\qquad b=1,\ldots,kq,\, (b,kq)=1,
\end{equation}
where
$$E_b(\sigma,(t_p)_{p>\exp{Q}})=Y_b^\varphi(\sigma,q)-\sum_{j=1}^d\sum_{h=2}^{\infty}\sum_{\substack{p> \exp{Q}\\p^h\equiv b\,(kq)}}\frac{f_j(p)^h}{hp^{h(\sigma+it_p)}},\qquad b=1,\ldots,kq, \,(b,kq)=1.$$
By \eqref{eq:ub_main_error} and \eqref{eq:ub_Y_b} we obtain
\begin{equation}\label{eq:ub_E_b}
|E_b(\sigma,(t_p)_{p>\exp{Q}})|\leq \frac{C_3}{q}, \quad b=1,\ldots,kq,\, (b,kq)=1,
\end{equation}
for some positive constant $C_3$, for every $Q^{-3/2}\leq \sigma-1\leq \eta$, if $r$, and thus $q$, is sufficiently large.\\
Let $R$ be the RHS of \eqref{eq:ub_E_b}. Then we look for solutions of the system
\begin{equation}\label{eq:linear_system}
\sum_{\substack{p> \exp{Q}\\p\equiv b\,(kq)}}\frac{f(p)}{p^{\sigma+it_p(\bb{z})}}=z_b,\quad b=1,\ldots,kq,\, (b,kq)=1,
\end{equation}
where $\bb{z}\in B_{R}(0)^{\phi(kq)}$ and $t_p:B_{R}(0)^{\phi(kq)}\rightarrow \R$ is a continuous function for every prime $p>\exp{Q}$, when $Q^{-3/2} \leq \sigma-1\leq \eta$ and $q$ is sufficiently large. Indeed, a solution of \eqref{eq:linear_system} gives a solution of \eqref{eq:sw_system} via Brouwer fixed point theorem in the following way (cf. Saias and Weingartner \cite[Lemma 2]{saias}). Consider the function
$$E: B_R(0)^{\phi(kq)}\rightarrow B_R(0)^{\phi(kq)}, \quad \bb{z}\mapsto (E_b(\sigma,(t_p(\bb{z}))_{p>\exp{Q}}))_{b=1,\ldots,kq,\,(b,kq)=1}.$$
By the absolute convergence of the Dirichlet series for $\sigma>1$ and the continuity of the functions $t_p(\bb{z})$ we have that $E(\bb{z})$ is continuous. Thus, by Brouwer fixed point theorem, there exists a fixed point $\bb{z}\in B_R(0)^{\phi(kq)}$, i.e.
$$E_b(\sigma,(t_p(\bb{z}))_{p>\exp{Q}})=z_b, \qquad b=1,\ldots, kq,\, (b,kq)=1.$$
Therefore, by \eqref{eq:linear_system}, we have found a solution of \eqref{eq:sw_system} for any $\sigma$ for which \eqref{eq:linear_system} holds.

To find solutions for \eqref{eq:linear_system} we observe that (cf. Jessen and Wintner \cite[Theorem 9]{jessenwintner})
$$\left\{\sum_{\substack{p> \exp{Q}\\p\equiv b\,(kq)}}\frac{f(p)}{p^{\sigma+it_p}}\,:\, t_p\in\R\right\}=\left\{z\,:\, \max\!\left(0,\max_{\substack{p> \exp{Q}\\p\equiv b\,(kq)}}\frac{2|f(p)|}{p^\sigma}-\sum_{\substack{p> \exp{Q}\\p\equiv b\,(kq)}}\frac{|f(p)|}{p^{\sigma}}\right)\leq |z|\leq\sum_{\substack{p> \exp{Q}\\p\equiv b\,(kq)}}\frac{|f(p)|}{p^{\sigma}}\right\}.$$
Thus the above plane set would be a disk if we had
\begin{equation*}\label{eq:circle}
\sum_{\substack{p> \exp{Q}\\p\equiv b\,(kq)}}\frac{|f(p)|}{p^{\sigma}}\geq \frac{2d}{\exp{Q\sigma}}, \qquad b=1,\ldots, kq, \, (b,kq)=1.
\end{equation*}
Suppose now that, if $q$ is sufficiently large, we can find $\sigma>1$ such that $Q^{-3/2} \leq \sigma-1\leq \eta$ and
\begin{equation}\label{eq:big_enough}
\sum_{\substack{p> \exp{Q}\\p\equiv b\,(kq)}}\frac{|f(p)|}{p^{\sigma}}\geq \max\left(18R,\frac{2d}{\exp{Q\sigma}}\right),\qquad \hbox{for every }b=1,\ldots,kq,\,(b,kq)=1.
\end{equation}
Then, for any $b=1,\ldots,kq$, $(b,kq)=1$, we take $p_{2,b}>p_{1,b}>\exp{Q}$ such that 
$$\sum_{\substack{\exp{Q}<p<p_{1,b}\\p\equiv b\,(kq)}}\frac{|f(p)|}{p^{\sigma}}<\frac{1}{3}\sum_{\substack{p> \exp{Q}\\p\equiv b\,(kq)}}\frac{|f(p)|}{p^{\sigma}} \leq \sum_{\substack{\exp{Q}<p\leq p_{1,b}\\p\equiv b\,(kq)}}\frac{|f(p)|}{p^{\sigma}},$$
and
$$\sum_{\substack{p_{1,b}<p<p_{2,b}\\p\equiv b\,(kq)}}\frac{|f(p)|}{p^{\sigma}}<\frac{1}{3}\sum_{\substack{p> \exp{Q}\\p\equiv b\,(kq)}}\frac{|f(p)|}{p^{\sigma}} \leq \sum_{\substack{p_{1,b}<p\leq p_{2,b}\\p\equiv b\,(kq)}}\frac{|f(p)|}{p^{\sigma}}.$$
We hence define
$$\mu_{1,b} = \frac{\sum_{\substack{\exp{Q}<p<p_{1,b}\\p\equiv b\,(kq)}}\frac{|f(p)|}{p^{\sigma}}}{\sum_{\substack{p> \exp{Q}\\p\equiv b\,(kq)}}\frac{|f(p)|}{p^{\sigma}}},\quad \mu_{2,b}=\frac{\sum_{\substack{p_{1,b}<p<p_{2,b}\\p\equiv b\,(kq)}}\frac{|f(p)|}{p^{\sigma}}}{\sum_{\substack{p> \exp{Q}\\p\equiv b\,(kq)}}\frac{|f(p)|}{p^{\sigma}}}\quad\hbox{and}\quad\mu_{0,b} = \frac{\frac{|f(p_{1,b})|}{p_{1,b}^\sigma}+\sum_{\substack{p\geq p_{2,b}\\p\equiv b\,(kq)}}\frac{|f(p)|}{p^{\sigma}}}{\sum_{\substack{p> \exp{Q}\\p\equiv b\,(kq)}}\frac{|f(p)|}{p^{\sigma}}}.$$
Then we have that the maps $G_b:(0,\pi/2)^2\rightarrow \C$, $(\theta_1,\theta_2)\mapsto\mu_{1,b}\exp{i\theta_1}+\mu_{2,b}\exp{-i\theta_2}$ are diffeomorphisms onto their image. Moreover, since $\mu_{1,b}+\mu_{2,b}-\mu_{0,b}>1/9$ and $|\mu_{1,b}-\mu_{2,b}|<2\exp{-Q}$, we have that (cf. Figure 1 of Saias and Weingartner \cite{saias})
$$\im{G_b}\supset \{w_b\in\C\, :\, |w_b-\mu_{0,b}|\leq 1/18 \}.$$
By \eqref{eq:big_enough} we may take
$$w_b=\mu_{0,b}+\frac{z_b}{\sum_{\substack{p> \exp{Q}\\p\equiv b\,(kq)}}\frac{|f(p)|}{p^{\sigma}}}$$
for any $|z_b|\leq R$, i.e. we may find a continuous solution
$$t_p(\bb{z})=\left\{\begin{array}{ll}
-\theta_1(z_b)/\log p & \exp{Q}<p<p_{1,b},\, p\equiv b\,(kq)\\
\theta_2(z_b)/\log p & p_{1,b}<p<p_{2,b},\, p\equiv b\,(kq)\\
\pi/\log p & p>p_{2,b}\hbox{ or }p=p_{1,b},\, p\equiv b\,(kq)\\
\end{array}\right.$$
of \eqref{eq:linear_system} for any $Q^{-3/2}\leq \sigma-1\leq \eta$ so that \eqref{eq:big_enough} holds, if $q$ is large enough. Therefore we just need to show that such $\sigma$ exists when $r$, and thus $q$, is large enough. First, we note that the second inequality in \eqref{eq:big_enough} trivially holds for any $\sigma\geq 1$ if $r$, and thus $q$, is sufficiently large. For the first inequality, by \eqref{eq:lb_main_term} and \eqref{eq:ub_E_b}, it is enough to find $1+Q^{-3/2}\leq \sigma\leq 1+\eta$ such that
\begin{equation*}\label{eq:last_ineq}
-\frac{C_1}{Q}+C_2\int_{ \exp{Q}}^\infty \frac{dx}{x^\sigma\log x}\geq 18C_3,
\end{equation*}
if $r$, and thus $q$, is sufficiently large. Using the asymptotic expansion of the exponential integral function (see e.g. Abramowitz and Stegun \cite[5.1.10]{abramowitzstegun}) we have
$$\int_{ \exp{Q}}^\infty \frac{dx}{x^\sigma\log x} = -\gamma - \log(Q(\sigma-1))-\sum_{k=1}^\infty \frac{Q^k(1-\sigma)^k}{k\cdot k!}.$$
It is then easy to check that the above inequality is verified if
$$\sigma-1\ll \frac{1}{Q}$$
and $r$, and thus $q$, is sufficiently large. We hence conclude that there exists $\sigma$ which verifies \eqref{eq:big_enough} if $r$, and thus $q$, is sufficiently large.\hfill\qed

%% file: continuity.tex
\section{Proof of Theorem \ref{theorem:continuity}}

We first note that if $\delta>0$ and $|\lambda-x_m|<\delta^2$ we have
\begin{equation*}
\abs{\sum_{n\leq 1/\delta} \frac{f(n)[\e{\lambda n}-\e{x_m n}]}{n^s} }\leq 2\pi \delta \sum_{n\leq 1/\delta} \frac{|f(n)|}{n^\sigma}.
\end{equation*}
Hence, by the absolute convergence of the Dirichlet series $F(s)$, for any $\sigma_0>1$ and any $\eta>0$ there exists $\delta_0=\delta_0(\eta,\sigma_0)$ such that
\begin{equation}\label{eq:num}
\abs{\sum_{n\leq 1/\delta} \frac{f(n)[\e{\lambda n}-\e{x_m n}]}{n^s} }< \frac{\eta}{4}, \qquad \hbox{for }\sigma>\sigma_0,\, \delta>\delta_0,\, |\lambda-x_m|<\delta^2.
\end{equation}

We now fix arbitrarily $\eps>0$ and $1<\sigma_0<\sigma^*+\eps$, where for simplicity $\sigma^*=\sigma^*(F,\lambda)$. Note that $\sigma^*>1$ by Theorem \ref{theorem:zero_lin_twist} when $\lambda$ is irrational while by Theorem 3 of \cite{righetti} when $\lambda=a/q$ is rational, since the family of functions $F(s,\chi)$, $\chi$ Dirichlet character mod $q$, satisfies the properties (E1)--(E5) of \cite{righetti} by the hypotheses on $F(s)$. Since $F(\lambda,1,1,s)\rightarrow f(1)\e{\lambda}\neq 0$ for $\sigma\rightarrow\infty$, by Bohr almost periodicity there exists $\eta>0$ such that $|F(\lambda,1,1,s)|\geq \eta$ for $\sigma\geq \sigma^*+\eps$ (cf. Bohr \cite[\S105]{bohr6}). Furthermore, by the the uniform convergence of the Dirichlet series $F(\lambda,1,1,s)$ for $\sigma\geq \sigma_0$, there exists $\delta_1=\delta_1(\eta,\sigma_0)>\delta_0$ such that for any $\delta>\delta_1$ we have
\begin{equation*}
\abs{\sum_{n\leq 1/\delta} \frac{f(n)\e{\lambda n}}{n^s} }\geq \frac{\eta}{2}, \qquad \hbox{for }\sigma\geq\sigma^*+\eps.
\end{equation*}
Then, by \eqref{eq:num} and the triangle inequality, for any $\delta>\delta_1$ and any $x_m$ such that $|\lambda-x_m|\leq \delta^2$ we get
\begin{equation*}
\abs{\sum_{n\leq 1/\delta} \frac{f(n)\e{x_m n}}{n^s} }\geq \frac{\eta}{4}, \qquad \hbox{for }\sigma\geq\sigma^*+\eps.
\end{equation*}
Finally, by the the uniform convergence of the Dirichlet series $F(x_m,1,1,s)$ for $\sigma\geq \sigma_0$, there exists $\delta_2=\delta_2(\eta,\sigma_0)>\delta_1$ such that for any $x_m$ with $|\lambda-x_m|<\delta_2^2$ we have $|F(x_m,1,1,s)|\geq \eta/8$ for $\sigma\geq\sigma^*+\eps$. Since $x_m\rightarrow \lambda$ as $m\rightarrow\infty$, we got that there exists $M_0$ such that for any $m\geq M_0$ we have $\sigma^*(F,x_m)<\sigma^*+\eps$.

Since $\sigma^*>1$ we may suppose that $\eps<\sigma^*-1$ and we may take $\rho=\beta+i\gamma$ such that $F(\lambda,1,1,\rho)=0$ with $\sigma^*-\eps<\beta\leq \sigma^*$. Furthermore we may suppose that $\sigma_0<\sigma^*-\eps$ and we fix $r>0$ such that $\sigma^*-\eps<\beta-r$ and $F(\lambda,1,1,s)\neq 0$ for $|s-\rho|=r$. Then, by \eqref{eq:num} and the uniform convergence of  $F(\lambda,1,1,s)$ and $F(x_m,1,1,s)$, there exists $\delta_3>0$ such that for every $x_m$ with $|\lambda-x_m|<\delta_3^2$ we have
\begin{spliteq*}
|F(\lambda,1,1,s)-F(x_m,1,1,s)|&\leq \abs{F(\lambda,1,1,s)-\sum_{n\leq 1/\delta_3}\frac{f(n)\e{\lambda n}}{n^s} }+\abs{\sum_{n\leq 1/\delta_3}\frac{f(n)\e{\lambda n}}{n^\sigma}-\frac{f(n)\e{x_m n}}{n^s} }\\
&\qquad \qquad+\abs{\sum_{n\leq 1/\delta_3}\frac{f(n)\e{\lambda n}}{n^s}-F(x_m,1,1,s) }\\
&< \min_{|s-\rho|=r}|F(\lambda,1,1,s)|
\end{spliteq*}
for every $\sigma\geq \sigma_0$, so in particular for $|s-\rho|\leq r$. Therefore, by Rouch\'e's theorem and since $x_m\rightarrow\lambda$ we get that there exists $M_1\geq M_0$ such that for any $m\geq M_1$ there exists $\rho_m=\beta_m+i\gamma_m$ with $F(x_m,1,1,\rho_m)=0$ and $\beta_m\geq \beta-r$. So, in particular, we got that $\sigma^*(F,x_m)\geq \beta-r>\sigma^*-\eps$ for every $m\geq M_1$.\hfill\qed